\documentclass[11pt]{amsart}

\usepackage{amssymb}
\usepackage[T1]{fontenc}
\usepackage{lmodern,microtype}
\input{glyphtounicode}
\usepackage{graphicx,tikz,placeins,needspace}
\usetikzlibrary{arrows.meta,positioning,calc}
\usepackage[hidelinks]{hyperref}
\newtheorem{theorem}{Theorem}[section]
\newtheorem{proposition}[theorem]{Proposition}
\newtheorem{lemma}[theorem]{Lemma}

\newtheorem{conjecture}[theorem]{Conjecture}
\theoremstyle{definition}

\theoremstyle{remark}
\newtheorem{remark}[theorem]{Remark}
\newcommand{\R}{\mathbb R}
\newcommand{\Z}{\mathbb Z}
\newcommand{\T}{\mathbb T}
\newcommand{\E}{\mathcal E}
\newcommand{\eps}{\varepsilon}
\newcommand{\dd}{\,d}
\newcommand{\ee}{\mathrm e}

\numberwithin{equation}{section}
\allowdisplaybreaks

\hypersetup{pdftitle={On the Bourgain-Rudnick eigenfunction restriction conjecture in four dimensions},pdfauthor={Ciprian Demeter}}
\title[Bourgain--Rudnick eigenfunction estimate in 4 dimensions]
{On the Bourgain--Rudnick eigenfunction restriction conjecture in four dimensions}
\author{Ciprian Demeter}
\address{Department of Mathematics, Indiana University,
	831 East 3rd St., Bloomington, IN 47405}
\email{demeterc@iu.edu}
\thanks{The author is partially supported by NSF grant DMS-2349828.}
\date{}

\begin{document}
	\begin{abstract}
		We prove a uniform upper $L^2$ bound for eigenfunctions of
		the four-dimensional flat torus on the compact paraboloid $\mathbb P^3$.
		The main novelty is the use of Fourier restriction estimates.
		We interpret spherical lattice point statistics as 
		non-concentration inequalities of the strength needed to apply  the Du--Zhang fractal restriction
		theorem. We apply their result on two-dimensional slices of $\mathbb P^3$. Further possible extensions are also discussed. 
	\end{abstract}
	\maketitle

	\section{Introduction and proof overview}\label{sec:intro}
	Bourgain and Rudnick formulated the following $L^2$ restriction inequality 
	\cite[Conjecture 1.9, (1.26)]{BR}. A reverse inequality is also postulated there, under curvature assumptions on $\Sigma$. That inequality will not concern us here. Write, for $R\ge1$ and $n\ge2$,
	$$\E_{R}^{(n)}=\Z^n\cap RS^{n-1},\qquad \ee(t)=e^{2\pi i t}.$$
	\begin{conjecture}[Bourgain--Rudnick]\label{conj:BR}
		For every $n\ge2$ and every fixed compact real analytic hypersurface
		$\Sigma\subset\T^n$ with surface measure $d\sigma_\Sigma$, there is a constant $C_\Sigma$ such that every
		toral Laplace eigenfunction 
		$$u(X)=\sum_{k\in\E_{R}^{(n)}}c_k\ee(k\cdot X)$$
		satisfies
		\begin{equation}\label{eq:BR-upper}
			\|u\|_{L^2(\Sigma)}\le C_\Sigma\|u\|_{L^2(\T^n)}.
		\end{equation}
		The constant is independent of $R$.
	\end{conjecture}
	No curvature assumption is made on $\Sigma$.
	Dependence on $\Sigma$ is unavoidable, as applying \eqref{eq:BR-upper} to a single Fourier mode gives
	$C_\Sigma\ge(\sigma_\Sigma(\Sigma))^{1/2}$.

	In \cite{BR} Bourgain--Rudnick established  \eqref{eq:BR-upper} in dimensions
	two and three for real analytic hypersurfaces with nonzero Gaussian
	curvature (equivalently, a nondegenerate second fundamental form). They also remark that in the absence of curvature, even when $n=2$ and $\Sigma$ is a geodesic segment, the conjectured inequality poses significant additional difficulties related to the distribution of lattice points in small caps. On the other hand, curvature allows one to use weaker arithmetic input. They exploit curvature by means of the detailed oscillatory structure of the stationary-phase expansion of $\widehat{d\sigma_\Sigma}$.    
	
	In dimension three, Bourgain--Rudnick 
	apply the $T^*T$ method to the restriction operator. Stationary phase
	expresses the resulting matrix entries
	$\widehat{d\sigma_\Sigma}(k-l)$ as oscillatory terms with
	$|k-l|^{-1}$ decay, up to smaller errors. Arithmetic cap estimates
	control nearby frequency pairs and these errors, while bounds for
	double exponential sums exploit cancellation among more distant pairs.

	The equivalent dual formulation of \eqref{eq:BR-upper} is
	\begin{equation}
		\label{eq:dual}
		\left(\sum_{k\in\E_R^{(n)}}
		|\widehat{g\,d\sigma_\Sigma}(k)|^2\right)^{1/2}
		\lesssim_\Sigma\|g\|_{L^2(\Sigma,d\sigma_\Sigma)}.
	\end{equation}
	Here $g\in L^2(\Sigma,d\sigma_\Sigma)$ and
	$\widehat{g\,d\sigma_\Sigma}(\xi)
	=\int_\Sigma g(X)\ee(-\xi\cdot X)\,d\sigma_\Sigma(X)$.

	Their Section~11 isolates the statistics
	used for frequency points. Their argument proves \eqref{eq:dual} for arbitrary point sets $\Lambda_R$ on $RS^2$ in place of $\E_R^{(3)}$, satisfying these statistics. Integer coordinates enter only in verifying the required statistics. In particular, their result is applicable to the idealized model of $\E_R^{(3)}$ consisting of $R^{1/2}$-separated points on $RS^2$, as the needed (in fact much stronger) upper bounds per cap are trivially satisfied in this case. The separation $R^{1/2}$ permits $\sim R$ points, matching the
	expected order of magnitude of $\#\E_R^{(3)}$ at typical radii. In general, the idealized model separation on $RS^{n-1}$ is $R^{\frac1{n-1}}$. In Section~11 they show that \eqref{eq:dual} fails for this model when $n\ge 8$, even when $\Sigma$ itself is $S^{n-1}$.

	\subsection{The new theorem}
	Let $\T^4=\R^4/\Z^4$ have normalized Haar measure, and drop the upper index to write
	\[
	\E_R=\Z^4\cap RS^3.
	\]
	An eigenfunction of $-\Delta$ with eigenvalue $4\pi^2R^2$ has the form
	\begin{equation}\label{eq:eigenfunction}
		u(X)=\sum_{k\in\E_R}c_k\ee(k\cdot X),\qquad
		\|u\|_{L^2(\T^4)}^2=\sum_{k\in\E_R}|c_k|^2.
	\end{equation}
	We use its periodic lift on the physical paraboloid
	\[
	\mathbb P^3=\{X(x):x\in [-1,1]^3\},\qquad X(x)=(x,|x|^2/2).
	\]
	We will prove the following result.
	
	\begin{theorem}\label{thm:main}
		There is a constant $C$, independent of $R$ and of the coefficients,
		such that
		\begin{equation}\label{eq:main}
			\left(\int_{[-1,1]^3}|u(X(x))|^2\dd x\right)^{1/2}
			\le C\|u\|_{L^2(\T^4)}.
		\end{equation}
	\end{theorem}

	Theorem~\ref{thm:main} gives the upper bound
	\eqref{eq:BR-upper} in dimension four for $\mathbb P^3$.
	Indeed, its surface measure is
	$d\sigma_{\mathbb P^3}=\sqrt{1+|x|^2}\,dx$, whose density lies between
	$1$ and $2$ on $[-1,1]^3$.

	The Bourgain--Rudnick analytic toolbox for $n=3$ combines the stationary-phase expansion of $\widehat{d\sigma_\Sigma}$ with classical $L^2$ methods, including
	Schur's test and rather fascinating Plancherel-based estimates for oscillatory sums.
	
	The major methodological difference in our argument is our use of band operators $T_a$ and of modern
	\emph{weighted Fourier restriction estimates} to obtain
	a saving for each localized operator. More precisely, we use as a black box the fractal $L^2$ restriction estimate of Du--Zhang \cite{DZ} in two frequency variables, applied to two-dimensional paraboloid slices of $\mathbb{P}^3$. The proof of their theorem packages the
	$\ell^2$ decoupling theorem  \cite{BDdec},
	broad--narrow analysis, and multilinear refined Strichartz estimates
	based on multilinear Kakeya. What allows for this application is our interpretation of lattice cap statistics (on the frequency side) as multi-scale non-concentration inequalities on the spatial side, as needed in the Du--Zhang theorem. This is achieved via duality. In contrast to Bourgain--Rudnick, we never explicitly use either the decay or the oscillatory profile of $\widehat{d\sigma_\Sigma}$.
	
	As in the Bourgain--Rudnick approach, our analytic argument does not require lattice structure. It works for any finite
	$\Lambda_R\subset RS^3$ satisfying \eqref{eq:intro-counts}, with
	$\E_R$ replaced by $\Lambda_R$. Integer coordinates enter in proving
	those counts and in the toral interpretation. The complexity of \eqref{eq:intro-counts} is roughly comparable to the estimates used by Bourgain--Rudnick.

	We present the argument for $\Sigma=\mathbb{P}^3$. Rather conveniently, its vertical
	slices are two-dimensional paraboloids, so the Du--Zhang theorem applies
	directly. See \eqref{eq:intro-shift}. For $\Sigma=S^3$, the corresponding slices are two-dimensional
	spheres, for which Kim's extension \cite{Kim} supplies the analogous
	weighted estimate.
	The same number-theoretic input and the overall approach of our paper
	suggest a route to smooth hypersurfaces with positive definite second
	fundamental form. Kim's weighted decoupling theorem \cite{Kim}
	generalizes the Du--Zhang theorem to strictly convex hypersurfaces. However, slicing a generic definite hypersurface introduces
	additional phases that depend nonlinearly on the projected frequencies
	and vary along the slices. The resulting operators are variable-phase
	oscillatory integral operators of Fourier-integral type, to which Kim's
	fixed-surface estimates do not apply directly. One would need to adapt
	the underlying decoupling and rescaling arguments to these phases while
	preserving the dependence on the two lattice cap statistics. Such a
	modification may address the positive definite case using the arithmetic
	estimates and the localization and summation strategy developed here.

	For a general nondegenerate hypersurface, the second fundamental form
	is invertible but can have mixed signs. The projected slices can then
	have indefinite or even degenerate curvature, so the elliptic input
	from Kim's theorem is no longer available. This prevents an immediate
	extension of the present argument to that setting.
	
	The zero curvature case looks to be even more challenging. For comparison, Zhang--Zhu \cite[Theorem~1.3(i)]{ZZ} obtain
	$R^{1/8+\eps}$ restriction bounds for totally geodesic hypersurface
	patches in $\T^4$, uniformly in their direction.

	\subsection{The main steps of the proof}
	Fix a radial $\eta\in C_c^\infty(B(0,2))$ with $0\le\eta\le1$,
	equal to one on $[-1,1]^3$, and define
	\begin{equation}\label{eq:T}
		T:\ell^2(\E_R)\longrightarrow L^2(\R^3,dx),\qquad
		Tc(x)=\eta(x)\sum_{k\in\E_R}c_k\ee(k\cdot X(x)).
	\end{equation}
	The affine normal $N(x)=(x,-1)$ will be convenient for the band
	decomposition. Here are the main steps of our argument.
	\medskip
	
	\emph{Upper bounds for lattice points on caps.}
	Section~\ref{sec:arithmetic} proves, for $1\le D\le R^{2/3}$,
	\begin{equation}\label{eq:intro-counts}
		\begin{aligned}
			\#(\E_R\cap B(z,r))&\lesssim_\eps D^\eps r^{7/4}
			\quad(1\le r\le CD),\\
			\#(\E_R\cap B(z,C\sqrt D))&\lesssim_\eps D^{1/2+\eps}.
		\end{aligned}
	\end{equation}
	In particular, the idealized model with $R^{1/3}$-separated points on $RS^3$ is easily seen to satisfy these statistics.
	The proof uses quadratic congruences, a determinant bound, the divisor estimate, and Minkowski's second theorem. 
	
	\emph{The main analytic input.} Section~\ref{sec:sieve} states Du--Zhang's result and proves its
	discrete restriction consequence, Proposition \ref{prop:slice}.

	\emph{Normal bands and exact slices.}
	Section~\ref{sec:local} decomposes $T$ into band operators
	\[
	T=T_0+\sum_{a\ge a_0}
	\bigl(T_{a,+}+T_{a,-}\bigr)
	\]
	$$
	T_ac(x)=\eta(x)\sum_{k\in\E_R}c_k\ee(k\cdot X(x))
	\chi\!\left(\frac{N(x)\cdot k}{Ra}\right)
	$$
	down to $a_0=R^{-1/3}$.  We further decompose $T_a$ into localized operators
	$$T_a=\sum_jT_aP_{Q_j},$$
	where $Q_j$ are cubes partitioning $\R^4$, with side length $D=a^{-2}$.
	
	Here is the effect of $P_Q$. If $(T_aP_Qc)(s,y)\ne0$ for some
	coefficient family $c$, then, after a horizontal rotation depending on $Q$,
	the coordinate $s$ lies in a common interval $I$ of length $O(a)$.
	We estimate each slice uniformly and then integrate over $I$. No
	cancellation in $s$ is needed. Here are some details.
	
	Fourier inversion expresses the band cutoff as an average of phase shifts. We then use the identity 
	\begin{equation}\label{eq:intro-shift}
		X(x)+hN(x)=((1+h)x,|x|^2/2-h).
	\end{equation}
	Freezing the $s$ coordinate of $x$ leaves a two-dimensional paraboloid. This is the part of the proof where 
	working with $\Sigma=\mathbb{P}^3$ makes things easier. 
	For each such slice, we apply the discrete fractal estimate in Proposition \ref{prop:slice}. The frequency points are projections to $\R^3$ of lattice points lying on a fixed cap. \eqref{eq:intro-counts} provides the  non-concentration hypotheses in the discrete restriction result.  Using Fubini in $s$ and the slice estimates leads to a good upper bound for the operator norm of $T_aP_Q$.

	\emph{Schur's test and final summation.}
	Section~\ref{sec:kernel} estimates the norm of $T_a$, via $\|T_a\|=\|T_a^*T_a\|^{1/2}$, the decomposition 
	$$T_a^*T_a=\sum_{i,j}(T_aP_{Q_i})^*(T_aP_{Q_j})
	=\sum_{i,j}K_{ij},$$
	the two bounds          
	\[
	\|K_{ij}\|\lesssim_L a^{1/16}\min\{1,a^{-6}(1+|i-j|)^{-L}\},
	\]
	and a final application of Schur's test. 
	The first bound follows from slicing and the Du--Zhang theorem. The second follows from kernel decay, proved by repeated integration by parts.

	\begin{remark}[$n=3$]
		In dimension three, normal bands can be used to give a simpler analytic proof of
		the Bourgain--Rudnick uniform upper bound \eqref{eq:BR-upper} for non-degenerate surfaces.
		No application of Du--Zhang is needed in this case.
	\end{remark}
	
	\begin{remark}[$n\ge 5$]
		The same analytic framework can also be used to prove \eqref{eq:BR-upper} in dimension five for $\mathbb P^4$, but requires
		substantially more involved number theory for the lattice set.
		The new feature is that this input forces the points in every sufficiently
		crowded cap on $RS^4$ to lie inside an affine hyperplane. After slicing,
		their contribution is controlled by a version of our $n=4$ result.
		Du--Zhang (this time for
		three-dimensional paraboloid slices) handles the remaining less crowded caps, as
		in our argument for $n=4$.
		
		This can be seen as an arithmetic analog of a broad--narrow argument,
		with the narrow part treated by dimension reduction, as in e.g.
		lower-dimensional decoupling. The additional structural input forcing
		lower-dimensional concentration relies critically on the lattice
		structure, as the idealized model lacks this rigidity. On the other hand,
		the idealized $R^{1/4}$-separated model in dimension five already has
		sufficiently small cap counts  and is also accessible to the
		analytic argument. This is consistent with the aforementioned
		construction of Bourgain and Rudnick that gives counterexamples
		for the idealized model when $n\ge8$.
		
		In dimension six, crowded five-dimensional sections and a diffuse
		contribution at critical cap occupancy create further obstacles.
		Our attempts have not resolved these difficulties, and we do not
		currently have a route through this method to the conjecture in
		general dimension.
	\end{remark}
	\subsection{Fourier decay in dimension three}\label{sec:n3-decay}
	The elementary three-dimensional analogue with an $R^\eps$ loss illustrates
	what Fourier decay alone can accomplish. If $\Sigma\subset\T^3$
	is a fixed compact smooth surface with nonzero Gaussian curvature,
	then every eigenfunction of frequency $R$ satisfies
	\begin{equation}\label{eq:n3-loss}
		\|u\|_{L^2(\Sigma)}
		\lesssim_{\Sigma,\eps}R^\eps\|u\|_{L^2(\T^3)}
		\qquad(\eps>0).
	\end{equation}
	Bourgain--Rudnick note this elementary bound in their introduction,
	before explaining the additional cancellation needed for uniformity
	\cite[Section 1, before Proposition 1.4]{BR}.

	We present this argument to explain the additional difficulty in
	dimension four and to motivate the band decomposition. In dimension three,
	the elementary estimate already reduces the problem to removing an
	arbitrarily small power loss. 
	
	Put $\E_R^{(3)}=\Z^3\cap RS^2$, and assume $R\ge2$.
	The only arithmetic input is
	\begin{equation}\label{eq:n3-count}
		\#(\E_R^{(3)}\cap B(z,r))\lesssim_\delta R^\delta r
		\qquad(r\ge1,\ \delta>0).
	\end{equation}
	Indeed, the first coordinate takes $O(r)$ integer values in the
	ball. For each value $k_1=m$, the other two coordinates satisfy
	$k_2^2+k_3^2=R^2-m^2$. For $1\le N\le R^2$, the
	sum-of-two-squares formula and the divisor bound give
	$r_2(N)\le4\tau(N)\lesssim_\delta R^\delta$; when $N=0$
	there is just one point. Summing over $m$ proves \eqref{eq:n3-count}.
	
	After a finite smooth localization, let $d\mu$ be a nonnegative
	smooth surface density on one patch of $\Sigma$. Stationary phase gives
	\[
	|\widehat\mu(\xi)|\lesssim_\Sigma(1+|\xi|)^{-1}.
	\]
	
	For $u(x)=\sum_{k\in\E_R^{(3)}}c_k\ee(k\cdot x)$, expanding
	the square gives the Gram form
	\[
	\int |u|^2\dd\mu
	=\sum_{k,l\in\E_R^{(3)}}c_k\overline{c_l}\,
	\widehat\mu(l-k).
	\]
	Since distinct lattice points have distance at least one,
	dyadic shells and \eqref{eq:n3-count} yield
	\[
	\begin{aligned}
		\sup_{k\in\E_R^{(3)}}\sum_{l\in\E_R^{(3)}}
		|\widehat\mu(l-k)|
		&\lesssim 1+
		\sup_k\sum_{j=0}^{\lceil\log_2(2R)\rceil}
		2^{-j}\#\bigl(\E_R^{(3)}\cap B(k,2^{j+1})\bigr)\\
		&\lesssim_\delta R^\delta\log(2R).
	\end{aligned}
	\]
	Symmetry gives the same column bound. Scalar Schur
	(Lemma~\ref{lem:Schur} below) therefore gives
	\[
	\int |u|^2\dd\mu
	\lesssim_{\Sigma,\delta}R^\delta\log(2R)\sum_k|c_k|^2.
	\]
	Choosing $\delta>0$ sufficiently small in terms of $\eps$, absorbing
	the logarithm, and taking square roots proves \eqref{eq:n3-loss}.

	In dimension four, the same absolute-value Gram--Schur argument produces
	an $R^{1/4}$ loss.
	
	\subsection{Precedents for the method}\label{sec:precedents}
	In their study of Fourier extension on large spheres and related oscillatory integral operators, Bennett--Seeger \cite[Sections~4--5]{BS} combine dyadic localization near a fold with lower-dimensional restriction estimates on input slices in model cases. This provides a methodological precedent for combining bands and slicing.

	Normal-component decompositions have precedents in semiclassical
	restriction. For semiclassical parameter $h$, Tacy
	\cite[Proposition~3.1]{Tacy} uses positive and negative normal bands
	together with a central $h^{1/3}$ region. See also Galkowski
	\cite[Section~3]{Galkowski}. Our explicit Fourier cutoffs combine
	this type of decomposition with arithmetic cap statistics and
	Du--Zhang estimates on exact slices.
	
	Arithmetic cap counting by well-chosen integer-normal slices
	already appears in Bourgain--Rudnick
	\cite[Lemma~2.3 and Section~2.6]{BR}.
	The successive-minima argument used below follows the same
	integer-form slicing principle as Zhang--Zhu
	\cite[Proposition~3.4]{ZZ}: select the forms with the smallest widths
	and fix their integer levels to obtain affine two-planes.
	Here we use the thin normal width of a spherical cap and the uniform
	local circle count, so the final loss depends only on the local radius.
	The material of
	Section~\ref{sec:circle-counting} follows from known results:
	Oravecz--Rudnick--Wigman \cite[(7.5)]{ORW} give the quadratic count,
	and Huang--Zhang \cite[Lemma~4]{HZ} give the circle count; three
	distinct lattice points within radius $r$ force their circle's radius to be
	$O(r^3)$. We include elementary proofs for completeness.

	Estimate~\eqref{eq:elementarycap} is also a short consequence of the
	affine-covering theorem of Balko--Cibulka--Valtr \cite[Theorem~5]{BCV},
	Minkowski's second theorem, and the uniform local circle bound.
	For a full-dimensional cap, the covering theorem supplies affine two-planes,
	Minkowski bounds their number, and the local circle bound keeps the loss in $r$, independently
	of $R$.
	We are not aware of its being recorded explicitly with a loss
	depending only on the local radius.

	\section{The arithmetic cap estimate}\label{sec:arithmetic}
	
	We prove the two counts needed by the weighted theorem. An elementary
	quadratic-congruence argument counts points on circles, and two
	integer linear forms with small widths on the cap slice a full-dimensional cap into circles. The
	thin normal width of the cap makes the number of slices small.
	
	\subsection{Circle counting with a loss in the local radius only}\label{sec:circle-counting}
	
	Oravecz--Rudnick--Wigman \cite[(7.5)]{ORW} give a stronger form
	of the following quadratic count. We include an elementary congruence proof.
	
	\begin{lemma}[Elementary quadratic count]\label{lem:normcount}
		For positive integers $d,N$,
		\[
		\#\{(X,Y)\in\Z^2:Y^2+dX^2=N\}\le16\tau(N)^2
		\lesssim_\delta N^\delta\qquad(\delta>0),
		\]
		uniformly in $d$. Here $\tau(N)$ is the number of positive divisors of $N$.
	\end{lemma}
	\begin{proof}
		Write $d=s^2d_0$ with $d_0$ squarefree. The injection
		$(X,Y)\mapsto(sX,Y)$ reduces to squarefree $d$.
		Writing $(X,Y)=g(x,y)$ with $\gcd(x,y)=1$ gives
		$g^2\mid N$ and $y^2+dx^2=n:=N/g^2$.
		For $n=1$ there are at most four solutions. For $n>1$,
		$\gcd(x,n)=1$, since a prime dividing $x$ and $n$ would divide $y$.
		Thus each primitive solution determines
		$t=yx^{-1}\pmod n$, with $t^2\equiv-d\pmod n$.
		
		There are at most $4\cdot2^{\omega(n)}$ possible residues $t$,
		where $\omega(n)$ counts the distinct prime divisors of $n$.
		If $p\mid\gcd(d,n)$, primitivity gives $p\nmid x$, $p\mid y$;
		squarefreeness forces $p$ to divide $n$ exactly once, with unique root $0$ modulo $p$.
		For odd $p\nmid d$,
		there are at most two roots modulo $p^h$: the difference of two
		root squares factors as $(t-u)(t+u)$, and at most one factor is
		divisible by $p$. For powers of $2$ with $d$ odd, there are at
		most four roots, because one of $t-u,t+u$ has exactly one factor
		of $2$. The Chinese remainder theorem gives the asserted bound.
		
		Each residue accounts for at most four primitive solutions.
		For two with the same residue, $xy'-x'y$ is a multiple of $n$.
		The vectors $(y,\sqrt d\,x)$ and $(y',\sqrt d\,x')$ both have
		length $\sqrt n$, so
		\[
		|xy'-x'y|\le n/\sqrt d.
		\]
		If $d>1$, this determinant is zero and the two vectors differ
		only by sign. If $d=1$, its possible values are $0,\pm n$;
		relative to a fixed vector, the other vector is either parallel
		or perpendicular, giving at most four choices. Therefore the
		primitive count is at most $16\cdot2^{\omega(n)}\le16\tau(n)$.
		Summing over $g^2\mid N$ gives $16\tau(N)^2$.
		The standard divisor estimate $\tau(N)\lesssim_\rho N^\rho$,
		with $2\rho\le\delta$, completes the proof.
	\end{proof}

	The following is a local form of Huang--Zhang
	\cite[Lemma~4]{HZ}, as explained in Section~\ref{sec:precedents}.
	We again include a direct proof.
	\begin{lemma}[Uniform local circle count]\label{lem:elementarycircle}
		Uniformly over circles $C$ in affine two-planes in $\R^4$ and
		balls $B$ of radius $r\ge1$,
		\begin{equation}\label{eq:circlecount}
			\#(C\cap B\cap\Z^4)\lesssim_\eps r^\eps.
		\end{equation}
	\end{lemma}
	\begin{proof}
		If there are at most two points, the assertion is immediate;
		this also handles points and line--sphere intersections.
		Otherwise choose three distinct integer points $p_0,p_0+u,p_0+v$
		in $C\cap B$. They are noncollinear, and $|u|,|v|\le2r$.
		After translating $p_0$ to zero, the relative circle center $w$
		lies in $\operatorname{span}(u,v)$ and satisfies
		$2w\cdot u=|u|^2$, $2w\cdot v=|v|^2$.
		Choose coordinates $i,j$ with
		$\Delta=u_i v_j-u_jv_i\ne0$, and put
		\[
		a=v_ju-u_jv,\qquad b=-v_iu+u_iv.
		\]
		An integer point $z$ of the translated circle is determined by
		$s=z_i,t=z_j$, since $z=(sa+tb)/\Delta$.
		The equation $|z|^2-2w\cdot z=0$ becomes
		\[
		As^2+2Bst+Ct^2+Ds+Et=0,
		\qquad A>0,\quad d:=AC-B^2>0.
		\]
		After multiplication by $\Delta^2$, these coefficients are
		$A=|a|^2$, $B=a\cdot b$, $C=|b|^2$,
		$D=-2\Delta w\cdot a$, and $E=-2\Delta w\cdot b$.
		They are integers of size $O(r^{C_0})$ for an absolute $C_0$: in particular,
		$2w\cdot a=v_j|u|^2-u_j|v|^2$, and similarly for $b$.
		Completing squares gives an injective map to
		\[
		\begin{gathered}
			X=2As+2Bt+D,\quad Y=2dt+AE-BD,\\
			Y^2+dX^2=N:=dD^2+(AE-BD)^2.
		\end{gathered}
		\]
		Both $d$ and $N$ have polynomial size in $r$. If $N=0$ there
		is at most one solution. Otherwise Lemma~\ref{lem:normcount} gives
		$O_\delta(N^\delta)$ possibilities. Choose $\delta$ small in terms
		of $\eps$ to obtain \eqref{eq:circlecount}.
	\end{proof}
	
	The point of anchoring at three nearby integer points is that
	the coefficients depend polynomially on $r$, even when the original
	circle has a very large radius or lies far from the origin.
	There is consequently no loss involving $R$ in this lemma.
	
	\subsection{Two slicing directions and two lattice counts}
	
	\begin{proposition}[Two lattice counts]\label{cor:density}
		Fix $C\ge1$. Uniformly for $1\le D\le R^{2/3}$,
		\begin{align}
			\#(\E_R\cap B(z,r))&\lesssim_{\eps,C}D^\eps r^{7/4},
			&&1\le r\le CD,\label{eq:density}\\
			\#(\E_R\cap B(z,C\sqrt D))&\lesssim_{\eps,C}
			D^{1/2+\eps}.\label{eq:occupancy}
		\end{align}
		The estimates hold for every center $z$ and for all subsets of the
		frequencies.
	\end{proposition}
	\begin{proof}
		We first prove the cap bound
		\begin{equation}\label{eq:elementarycap}
			\#(\E_R\cap B(z,r))\lesssim_\eps
			r^\eps\left(r+\frac{r^{5/2}}{\sqrt R}\right)
			\qquad(1\le r\le2R).
		\end{equation}
		
		Let $A=\E_R\cap B(z,r)$ and choose $k_0\in A$ if $A$ is nonempty.
		All differences $k-k_0$ have length at most $2r$.
		First suppose that $A$ has affine dimension four. With
		$\nu_0=k_0/R$, the identity
		\[
		\nu_0\cdot(k-k_0)=-\frac{|k-k_0|^2}{2R}
		\]
		places these differences in a symmetric box $K$ with three
		tangential widths $Cr$ and one normal width $Cr^2/R$.
		In particular,
		\begin{equation}\label{eq:elementaryvolume}
			V:=\operatorname{vol}_4(K)\asymp r^5/R.
		\end{equation}
		
		Define the width and polar body by
		\[
		w(\xi)=\sup_{z\in K}|\xi\cdot z|,
		\qquad K^\circ=\{\xi:w(\xi)\le1\}.
		\]
		The $j$th successive minimum is the least $t>0$ for which
		$tK^\circ$ contains $j$ linearly independent integer vectors.
		Minkowski's second theorem \cite[Chapter~VIII, Theorem~V]{Cassels}
		bounds the product of these minima
		by $2^4/\operatorname{vol}_4(K^\circ)$. Since
		$\operatorname{vol}_4(K^\circ)\asymp V^{-1}$ for a box,
		there are independent $\xi_j\in\Z^4$ with ordered widths
		\[
		w_1\le w_2\le w_3\le w_4,\qquad
		\prod_{j=1}^4w_j\lesssim V.
		\]
		Every nonzero integer vector has width at least one. Otherwise all
		its integer pairings with the spanning differences $k-k_0$ would
		vanish, contradicting their linear independence.
		Consequently,
		\begin{equation}\label{eq:elementarywidths}
			1\le w_1\le w_2,\qquad w_1w_2\lesssim V^{1/2}.
		\end{equation}
		
		Slice $A$ by the two integer levels
		\[
		\xi_1\cdot(k-k_0)=h_1,\qquad
		\xi_2\cdot(k-k_0)=h_2,\qquad h_1,h_2\in\Z.
		\]
		There are $O(w_1w_2)$ possible pairs of levels.
		Each common level is an affine two-plane; its intersection with
		$RS^3$ is a circle, a point, or empty. Apply
		Lemma~\ref{lem:elementarycircle} in the original radius-$r$ ball.
		The slicing vectors may have large coefficients; the circle lemma
		is independent of those coefficients, so this incurs no extra loss.
		Equations~\eqref{eq:elementaryvolume}--\eqref{eq:elementarywidths}
		give
		\[
		\#A\lesssim_\eps r^\eps V^{1/2}
		\lesssim_\eps r^\eps r^{5/2}R^{-1/2}.
		\]
		
		If $A$ has affine dimension three, choose a coordinate $k_j$
		which is nonconstant on its affine span. Since $A$ lies in a
		ball of radius $r$, the coordinate $k_j$ takes only $O(r)$
		distinct integer values on $A$. Each corresponding hyperplane
		$x_j=m$ intersects the affine span of $A$ in an affine two-plane.
		The circle lemma gives $\#A\lesssim_\eps r^{1+\eps}$.
		For affine dimension at most two, the circle lemma applies directly
		(or the sphere meets the affine span in at most two points).
		Combining the cases proves \eqref{eq:elementarycap}.

		Finally, $R\ge D^{3/2}$ gives
		\[
		\frac{r^{5/2}}{\sqrt R}\lesssim_C r^{7/4}
		\quad(1\le r\le CD),\qquad
		\frac{(C\sqrt D)^{5/2}}{\sqrt R}\lesssim_C D^{1/2}.
		\]
		Insert these in \eqref{eq:elementarycap} and absorb $r^\delta$
		into $D^\eps$ by choosing $\delta$ small. This proves both
		claims. Fixed radius enlargements and bounded $R$ only change
		the constants; counts decrease on taking subsets.
	\end{proof}

	\section{Du--Zhang and the discrete restriction estimate}\label{sec:sieve}
	
	We state the weighted extension theorem and pass from integrals
	to point samples by a short smoothing identity. Taking the adjoint
	gives a discrete estimate on one fixed paraboloid $\mathbb{P}^2$.
	
	\subsection{The weighted estimate}
	Define
	\[
	Z(y)=(y_1,y_2,|y|^2/2),\qquad
	Ef(p)=\int_{[-2,2]^2} f(y)\ee(-p\cdot Z(y))\dd y.
	\]
	
	Let $D\ge2$ and let $Y\subset B(0,CD)\subset\R^3$ be a union
	of distinct unit cubes from a fixed lattice grid. Suppose, for $1\le\alpha\le3$,
	\begin{equation}\label{eq:Ydensity}
		|Y\cap B(z,r)|\lesssim\gamma r^\alpha\quad(1\le r\le CD),
		\qquad |Y\cap B(z,C_1\sqrt D)|\lesssim\lambda.
	\end{equation}
	For nonempty $Y$ take $\gamma,\lambda\ge1$. Fixed radius
	enlargements only alter constants.
	
	\begin{theorem}[Du--Zhang, {\cite[Theorem~1.6]{DZ}}]\label{thm:DZ}
		Under \eqref{eq:Ydensity}, for every $\eps>0$,
		\begin{equation}\label{eq:weighted}
			\|Ef\|_{L^2(Y)}^2
			\lesssim_{\eps,\alpha,C,C_1}D^{\alpha/6+\eps}
			\gamma^{1/3}\lambda^{1/3}\|f\|_2^2.
		\end{equation}
	\end{theorem}
	
	\subsection{Sampling and duality}
	We will need the following discrete corollary.
	\begin{proposition}[Discrete restriction estimate]\label{prop:slice}
		Let $A\subset B(0,CD)\subset\R^3$ be finite. Suppose, for every
		$\delta>0$,
		\begin{equation}\label{eq:Adensity}
			\begin{aligned}
				\#(A\cap B(z,r))&\lesssim_\delta D^\delta r^{7/4}
				\quad(1\le r\le CD),\\
				\#(A\cap B(z,C_1\sqrt D))&\lesssim_\delta D^{1/2+\delta}.
			\end{aligned}
		\end{equation}
		Then
		\begin{equation}\label{eq:slice}
			\int_{[-2,2]^2}\left|\sum_{p\in A}d_p\ee(p\cdot Z(y))\right|^2\dd y
			\lesssim D^{15/32}\sum_{p\in A}|d_p|^2.
		\end{equation}
		The points of $A$ need not have integer coordinates.
	\end{proposition}
	
	Proposition \ref{prop:slice} reflects a reduction from $\R^4$ to $\R^3$ on both sides of the original problem. Spatially, slicing $\mathbb{P}^3$ by vertical hyperplanes gives two-dimensional paraboloids, and Fubini combines the slice estimates. In frequency, we project the spherical  lattice points in a small cap of $RS^3$ onto $\R^3$. These projections, up to controlled affine changes, form $A$.
	\begin{proof}
		Let $S_A:\ell^2(A)\to L^2([-2,2]^2,dy)$ denote
		the exponential-sum operator on the left. Then
		\[
		(S_A^*f)(p)=Ef(p),
		\]
		and $\|S_A\|=\|S_A^*\|$.
		Choose a nonnegative smooth $w$ on $\R^3$, of integral one,
		supported in a sufficiently small fixed ball. Its Fourier transform
		$\widehat w(z)=\int w(v)\ee(-v\cdot z)\dd v$ satisfies
		$|\widehat w(Z(y))|\ge1/2$ for $y\in[-2,2]^2$.
		For $g(y)=f(y)/\widehat w(Z(y))$, Fubini gives
		\[
		Ef(p)=\int w(v)Eg(p+v)\dd v,
		\qquad \|g\|_2\le2\|f\|_2.
		\]
		The radius-one case of \eqref{eq:Adensity} bounds the number of
		small balls about $A$ covering any point by $O_\delta(D^\delta)$.
		Cauchy--Schwarz in the preceding identity therefore gives
		\begin{equation}\label{eq:sampling}
			\sum_{p\in A}|Ef(p)|^2\lesssim_\delta
			D^\delta\int_Y|Eg(z)|^2\dd z,
		\end{equation}
		where $Y$ is the union of unit cubes meeting the small balls about
		$A$. Each point contributes only boundedly many cubes, so
		\eqref{eq:Adensity} gives
		\eqref{eq:Ydensity} with
		$\alpha=7/4$, $\gamma\lesssim_\delta D^\delta$, and
		$\lambda\lesssim_\delta D^{1/2+\delta}$.
		Du--Zhang, with loss $D^\delta$, bounds \eqref{eq:sampling} by
		\[
		D^\delta D^{7/24+\delta}D^{\delta/3}
		D^{1/6+\delta/3}\|f\|_2^2
		=D^{11/24+8\delta/3}\|f\|_2^2.
		\]
		Choose $\delta=1/256$; then $11/24+8\delta/3=15/32$.
		Taking the adjoint proves the assertion. This fixed exponent
		leaves the positive gap $1/2-15/32=1/32$ needed below.
	\end{proof}

	\section{Normal bands and exact paraboloid slices}\label{sec:local}
	
	We restrict the frequencies to a cap on the frequency sphere $RS^3$.
	The normal band is shown to confine one coordinate of the physical
	paraboloid $\mathbb P^3$ to an interval of length $O(a)$.
	After projection, the frequencies form a point set in $\R^3$ to
	which the discrete slice estimate applies.
	We may assume throughout the rest of the proof that $R$ exceeds a
	fixed sufficiently large constant. For bounded $R$, Cauchy--Schwarz and
	the bounded number of lattice frequencies give the theorem directly.
	
	\subsection{The bands}\label{sec:band-definition}
	Choose an even smooth $\psi_0$, equal to one on $[-1,1]$ and
	supported in $[-2,2]$. Set $\chi(t)=\psi_0(t/2)-\psi_0(t)$,
	and let $\chi_\pm$ be its restrictions to the positive and
	negative half-lines, extended by zero. Their supports lie in
	$[1,4]$ and $[-4,-1]$.

	Recall the normal 
	$N(x)=(x,-1)$ to $\mathbb{P}^3$.
	The parameter $a$ records the size of the
	normal component $N(x)\cdot k/R$. We use dyadic scales
	$a=2^ja_0$, $j\ge0$, where
	\begin{equation}\label{eq:scales}
		a_0=R^{-1/3}.
	\end{equation}
	A band operator is
	\begin{equation}\label{eq:band}
		T_ac(x)=\eta(x)\sum_{k\in\E_R}c_k\ee(k\cdot X(x))
		\chi\!\left(\frac{N(x)\cdot k}{Ra}\right).
	\end{equation}
	Here $\chi=\chi_+$ or $\chi_-$ gives $T_{a,+}$ or $T_{a,-}$;
	using $\psi_0$ at $a=a_0$ gives the central operator $T_0$.
	The exact decomposition is
	\[
	T=T_0+\sum_{j\ge0}
	\bigl(T_{2^ja_0,+}+T_{2^ja_0,-}\bigr).
	\]
	Indeed, the definition of $\chi$ telescopes, and
	$\psi_0(t/(2^ja_0))=1$ for all sufficiently large $j$.
	Only finitely many terms occur: on the support of $\eta$,
	$|N(x)\cdot k/R|\le\sqrt5<3$, so bands with $a>3$ vanish.
	
	The signed cutoffs fix both the sign
	and size of each frequency's normal component. The sign will rule out stationary
	phase for pairs of distinct frequencies in Section~\ref{sec:kernel}. 
	The size confines a cap's contribution to a slab of thickness $O(a)$. 
	
	We call $a_0\le a\le a_*:=1/100$ a \emph{small-band} scale and
	$a_*<a\le3$ a \emph{large-band} scale. The large bands are handled by
	Proposition~\ref{prop:kernel}. There are only $O(1)$ such dyadic scales.
	We now fix a small-band scale and decompose $T_a$ into localized pieces. 
	
	\subsection{The localized operators}
	For $Q\subset\R^4$, define
	\[
	P_Q:\ell^2(\E_R)\to\ell^2(\E_R),\qquad
	(P_Qc)_k=c_k\mathbf1_Q(k).
	\]
	Fix $a$ and let $D=a^{-2}$.
	Partition $\R^4$ into
	cubes $Q_i=D(i+[0,1)^4)$, $i\in\Z^4$.
	Every occupied cube lies in a ball $B(Re_i,4D)$, by choosing
	$Re_i\in\E_R\cap Q_i$.
	Write 
	$$T_a=\sum_jT_aP_{Q_j}.$$
	We will first estimate the norm of the local operators  $T_aP_{Q_j}$,  and then the norm of $T_a$. It is far from immediate that the norm of the local operators is even $\lesssim 1$. We will need something slightly better, and this will be proved in \eqref{eq:local}.

	Fix an occupied cube $Q$ and choose $Re\in\E_R\cap Q$, so that
	$Q\subset B(Re,4D)$. Note that
	$$T_aP_Qc(x)=\eta(x)\sum_{k\in\E_R\cap Q}c_k\ee(k\cdot X(x))
	\chi\!\left(\frac{N(x)\cdot k}{Ra}\right).$$
	\subsection{The common interval $I$}
	Write $e=(e',c)$ with $e'\in\R^3$, $c\in\R$, and put $b=|e'|$.
	If the summand indexed by $k\in\E_R\cap Q$ is nonzero at some $x$, then
	$|x|\le2$, $|N(x)\cdot k|\le4Ra$, and $|k-Re|\le4D$. Thus
	\begin{equation}\label{eq:active}
		\begin{aligned}
			|x\cdot e'-c|=|N(x)\cdot e|
			&\le |N(x)\cdot k|/R+|N(x)|\,|Re-k|/R\\
			&\le4a+4\sqrt5\,D/R\le16a,
			\qquad \text{ since }D\le Ra.
		\end{aligned}
	\end{equation}
	Recall that $a\ge R^{-1/3}$.
	Since $1\le b+|c|\le3b+16a$ and $16a\le 16a_*<1/2$, every $Q$
	whose contribution is nonzero satisfies $b\ge b_0:=1/6$.
	
	We henceforth consider such caps.
	Rotate the horizontal coordinates so that $e=(b,0,0,c)$ and
	write $x=(s,y)$, $y\in\R^2$. This preserves the paraboloid and
	the radial cutoff $\eta$. Every summand vanishes unless
	\begin{equation}\label{eq:interval}
		s\in I:=\left[\frac{c-16a}{b},\frac{c+16a}{b}\right],
		\qquad |I|=32a/b\lesssim a.
	\end{equation}
	This interval is independent of $k$, so the whole $Q$-cap contribution
	is supported in the common slab $I\times\R^2$. On the support of
	$\eta$, the remaining coordinates satisfy $y\in[-2,2]^2$.
	
	\begin{figure}[htb]
		\centering
		\resizebox{\linewidth}{!}{%
			\begin{tikzpicture}[scale=1.18,>=Stealth,line width=.65pt,
				draw=black,every node/.style={text=black,font=\normalsize}]
				\path[use as bounding box] (-1.95,-2.65) rectangle (9.0,2.8);
				
				\draw (0,0) circle[radius=1.65];
				\draw[dashed,dash pattern=on 4pt off 4pt]
				(-1.647,-.10) arc[start angle=180,end angle=0,x radius=1.647,y radius=.23];
				\draw (-1.647,-.10)
				arc[start angle=180,end angle=360,x radius=1.647,y radius=.23];
				\draw (.28,1.626)
				.. controls (.19,.93) and (.79,.40) .. (1.642,.162);
				\node at (-1.46,1.72) {$RS^3$};
				\node at (.68,1.04) {$Q$};
				\fill[black] (1.167,1.167) circle[radius=1.45pt];
				\node[anchor=west] at (1.36,1.49) {$Re=(Rb,0,0,Rc)$};
				\fill[black] (1.03,.57) circle[radius=1.35pt];
				\node[anchor=west,inner sep=2pt] at (1.13,.57) {$k$};
				
				\begin{scope}[xshift=6.5cm]
					\node at (.88,2.49) {$\mathbb P^3$};
					\draw plot[domain=-1.18:1.25,samples=80] (\x,{1.28*\x*\x});
					\draw[->] (0,0)--(2.13,0);
					\node[anchor=north west,inner sep=2pt] at (1.98,-.10) {$y$};
					\draw[->] (0,0)--(-1.15,-2.20);
					\node[anchor=north,inner sep=3pt] at (-1.15,-2.24) {$s$};
					
					\draw[dashed,dash pattern=on 4pt off 4pt]
					(-.382,-.73)--(1.65,-.73);
					\draw[dashed,dash pattern=on 4pt off 4pt]
					(-.585,-1.12)--(1.65,-1.12);
					\node[anchor=east,inner sep=2pt] at (-.64,-.92) {$I$};
					
					\draw[dashed,dash pattern=on 4pt off 4pt]
					(.96,1.179648)--(.96,-.925);
					\fill[black] (.96,1.179648) circle[radius=1.45pt];
					\draw[->] (.96,1.179648)--(1.94,.781);
					\node[anchor=south,inner sep=2pt] at (1.68,1.28) {$N(x)$};
					\fill[black] (.96,-.925) circle[radius=1.35pt];
					\node[anchor=west,inner sep=3pt] at (1.03,-.925) {$x$};
				\end{scope}
			\end{tikzpicture}%
		}
		\caption{The frequency cap $Q\cap RS^3$ and the physical paraboloid
			$\mathbb P^3$, with the common interval $I$ in the $s$ coordinate.}
		\label{fig:normalband}
	\end{figure}
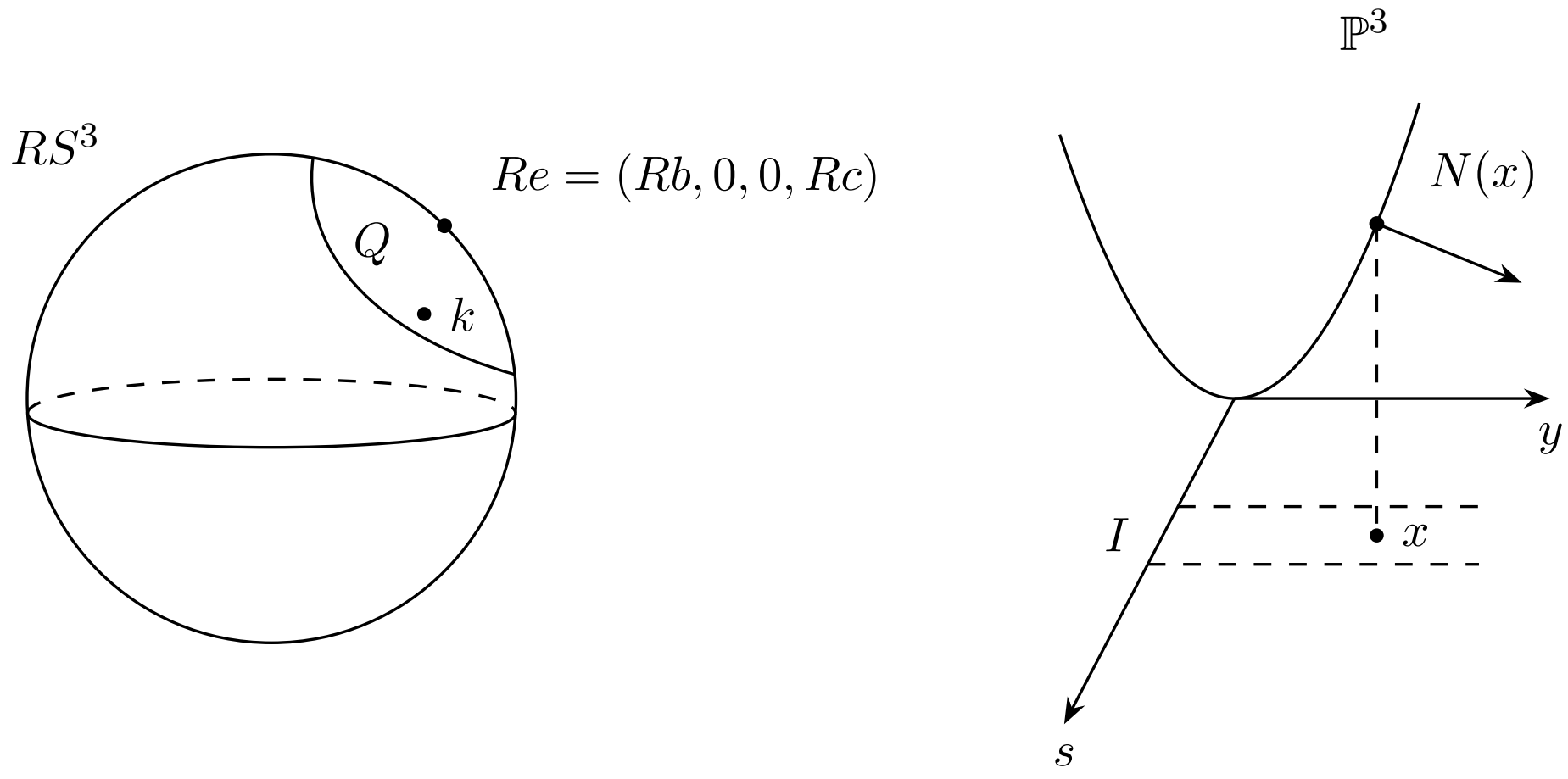

	\subsection{Projection and the reduction to the discrete slice estimate}
	In these rotated coordinates define the projection
	\[
	\pi:\R^4\longrightarrow\R^3,\qquad
	\pi(\xi_1,\xi_2,\xi_3,\xi_4)=(\xi_2,\xi_3,\xi_4).
	\]
	Since $Re=(Rb,0,0,Rc)$, its image is $p_0:=\pi(Re)=(0,0,Rc)$.
	For $k,l\in\E_R\cap Q$, write $k=(k_1,k_2,k_3,k_4)$ and
	$l=(l_1,l_2,l_3,l_4)$ in these coordinates: $k_1,l_1$ are precisely
	the coordinates dropped by $\pi$.
	For large $R$, since $D/R\le R^{-1/3}$,
	\[
	k_1,l_1\ge Rb-4D\ge Rb_0/2.
	\]
	Subtracting $k_1^2+|\pi k|^2=l_1^2+|\pi l|^2=R^2$ gives
	\[
	(k_1-l_1)(k_1+l_1)=-(\pi k-\pi l)\cdot(\pi k+\pi l).
	\]
	Consequently,
	\[
	|k_1-l_1|\le\frac{|\pi k|+|\pi l|}{k_1+l_1}
	|\pi k-\pi l|\le\frac2{b_0}|\pi k-\pi l|,
	\]
	and Pythagoras yields
	\[
	|\pi k-\pi l|\le|k-l|
	\le\sqrt{1+4/b_0^2}\,|\pi k-\pi l|\lesssim|\pi k-\pi l|.
	\]
	Thus projection is one-to-one and changes distances by fixed
	factors. The horizontal rotation need not preserve integer
	coordinates, but it preserves distances and hence all cap counts.
	For $|h|\le1/2$, the same is true of
	\begin{equation}\label{eq:projected}
		p_h(k)=((1+h)k_2,(1+h)k_3,k_4)-p_0,
		\quad A_h=\{p_h(k):k\in\E_R\cap Q\}.
	\end{equation}
	
	The sets $A_h$ lie in $B(0,CD)$.
	A ball of radius $r$ meeting $A_h$ lifts into a ball of radius
	$Cr$ about one of its preimages. Proposition~\ref{cor:density}
	therefore gives \eqref{eq:Adensity} for $A_h$, uniformly in $h$.

	\Needspace{8\baselineskip}
	For $F_Q(z)=\sum_{k\in\E_R\cap Q}c_k\ee(k\cdot z)$, the
	identity 
	$$
	X(x)+hN(x)=((1+h)x,|x|^2/2-h)
	$$
	gives
	for $x=(s,y)$
	\begin{equation}\label{eq:phase}
		k\cdot(X+hN)=(1+h)k_1s+k_4(s^2/2-h)
		+(p_h(k)+p_0)\cdot Z(y).
	\end{equation}
	These shifted evaluations are motivated by the Fourier inversion
	formula \eqref{eq:normal-average} below, which expresses the band
	operator as a weighted average of them, with $h=t/(Ra)$.
	At fixed $s,h$, put
	\[
	d_k=c_k\ee\bigl((1+h)k_1s+k_4(s^2/2-h)\bigr),
	\qquad |d_k|=|c_k|.
	\]
	The $p_0$ term gives a common factor of absolute value one.
	Proposition~\ref{prop:slice} therefore gives
	\[
	\begin{aligned}
		&\int_{[-2,2]^2}|F_Q(X(s,y)+hN(s,y))|^2\dd y\\
		&\quad=\int_{[-2,2]^2}
		\left|\sum_{k\in\E_R\cap Q}d_k\ee(p_h(k)\cdot Z(y))\right|^2\dd y\\
		&\quad\lesssim D^{15/32}\sum_{k\in\E_R\cap Q}|c_k|^2
		=D^{15/32}\|P_Qc\|_2^2.
	\end{aligned}
	\]
	Integrating over $s\in I$ and
	using $|I|\lesssim a=D^{-1/2}$ gives, for every $|h|\le1/2$,
	\begin{equation}\label{eq:uncut-local}
		\|\mathbf1_I(s)\eta(s,y)F_Q(X+hN)\|_{L^2(ds\,dy)}^2
		\lesssim D^{-\kappa}\|P_Qc\|_2^2,
		\quad \kappa=1/32.
	\end{equation}

	\subsection{The norm of $T_aP_Q$ for a small-band scale $a$}

	Fourier inversion writes the cutoff
	$\chi(N(x)\cdot k/(Ra))$ as an average of exponentials, allowing us to
	apply \eqref{eq:phase} and then \eqref{eq:uncut-local}. 
	\begin{proposition}[Local cap saving]\label{prop:local}
		For every occupied cube $Q\subset B(Re,4D)$ as above and every small signed or central band,
		\begin{equation}\label{eq:local}
			\|T_aP_Q\|_{\ell^2(\E_R)\to L^2(\R^3)}^2\lesssim D^{-\kappa},
			\qquad \kappa=1/32.
		\end{equation}
	\end{proposition}
	\begin{proof}
		With $\widehat\chi(t)=\int\chi(v)\ee(-tv)\dd v$, Fourier inversion gives
		
		\begin{equation}\label{eq:normal-average}
			\begin{aligned}
				T_aP_Qc(x)&=
				\eta(x)\sum_{k\in\E_R\cap Q}c_k\ee(k\cdot X(x))
				\chi\!\left(\frac{N(x)\cdot k}{Ra}\right)\\
				& =
				\eta(x)\int_{\R}\widehat\chi(t)
				F_Q\!\left(X(x)+\frac{t}{Ra}N(x)\right)\dd t.
			\end{aligned}
		\end{equation}
		
		Recall that we write $x=(s,y)$.
		The left side is supported in $I\times\R^2$, so insert
		$\mathbf1_I(s)$ on the right, keeping this factor in both pieces
		of the following split. For $|t|\le Ra/2$, the displacement
		$h=t/(Ra)$ satisfies $|h|\le1/2$. Minkowski and
		\eqref{eq:uncut-local} give
		\[
		\begin{aligned}
			&\left\|\mathbf1_I(s)\eta(x)\int_{|t|\le Ra/2}\widehat\chi(t)
			F_Q\!\left(X+\frac{t}{Ra}N\right)\dd t\right\|_{L^2(dx)}\\
			&\quad\le\int_{|t|\le Ra/2}|\widehat\chi(t)|
			\left\|\mathbf1_I\eta F_Q\!\left(X+\frac{t}{Ra}N\right)\right\|_{L^2(dx)}\dd t
			\lesssim D^{-\kappa/2}\|P_Qc\|_{2}\|\widehat\chi\|_1.
		\end{aligned}
		\]

		The remaining piece is the tail
		\[
		E_{a,Q}c(x):=\mathbf1_I(s)\eta(x)
		\int_{|t|>Ra/2}\widehat\chi(t)
		F_Q\!\left(X(x)+\frac{t}{Ra}N(x)\right)\dd t.
		\]
		Here  $|h|=|t|/(Ra)>1/2$, so the preceding slice
		estimate does not apply. Instead we use rapid decay of
		$\widehat\chi$ together with a bound for $F_Q$ valid at every
		point of $\R^4$. Cauchy--Schwarz in the frequency sum gives
		\[
		|F_Q(z)|\le\#(\E_R\cap Q)^{1/2}\|P_Qc\|_2
		\lesssim D^2\|P_Qc\|_2.
		\]
		Since $|I|\lesssim a$ and $\eta$ has fixed bounded support,
		$\|\mathbf1_I\eta\|_2\lesssim a^{1/2}$. Taking absolute values
		inside the tail integral therefore gives
		\[
		\begin{aligned}
			\|E_{a,Q}c\|_2
			&\le \|\mathbf1_I\eta\|_2\,
			\sup_{z\in\R^4}|F_Q(z)|
			\int_{|t|>Ra/2}|\widehat\chi(t)|\dd t\\
			&\lesssim a^{1/2}D^2\|P_Qc\|_2
			\int_{|t|>Ra/2}|\widehat\chi(t)|\dd t.
		\end{aligned}
		\]
		Smoothness and compact support of $\chi$ give rapid decay of
		$\widehat\chi$, hence
		$\int_{|t|>Ra/2}|\widehat\chi(t)|\dd t\lesssim_B(Ra)^{-B}$
		for every $B>0$. Substituting this bound gives
		\[
		\|E_{a,Q}c\|_2
		\lesssim_B a^{1/2}D^2(Ra)^{-B}\|P_Qc\|_2
		\le D^{7/4-B}\|P_Qc\|_2.
		\]
		Here $a^{1/2}=D^{-1/4}$ and $Ra\ge D$. Taking $B=2$ makes the
		tail norm $O(D^{-1/4})\|P_Qc\|_2$, smaller than the required
		$D^{-\kappa/2}\|P_Qc\|_2$. Combining the two pieces by the
		triangle inequality and then squaring proves the assertion.
		
	\end{proof}
	
	\section{Gram decay and Schur summation}\label{sec:kernel}\label{sec:global}
	
	The Gram matrix turns the problem of combining caps into an
	operator estimate on the frequency coefficients. We prove decay
	beyond distance $D=a^{-2}$, state and prove the Schur test, and
	combine this decay with the local saving. Summing the normal
	bands then completes the proof.
	
	\subsection{The Gram matrix controls restriction}\label{sec:gram}
	Write $T_ac=\sum_k c_kf_{a,k}$, where
	\[
	f_{a,k}(x)=\eta(x)\ee(k\cdot X(x))
	\chi\!\left(\frac{N(x)\cdot k}{Ra}\right).
	\]
	With inner products linear in the first argument, the adjoint is
	$(T_a^*g)_k=\langle g,f_{a,k}\rangle$.
	Thus $K_a=T_a^*T_a$ is a linear operator on $\ell^2(\E_R)$ associated with the \emph{Gram matrix} of these functions. 
	Its entries are their pairwise inner products,
	\begin{equation}\label{eq:kerneldef}
		K_a(k,l)=\langle f_{a,l},f_{a,k}\rangle
		=\int_{\R^3}A_{a,k,l}(x)\ee((l-k)\cdot X(x))\dd x,
	\end{equation}
	where
	\[
	A_{a,k,l}(x)=|\eta(x)|^2
	\overline{\chi\!\left(\frac{N(x)\cdot k}{Ra}\right)}
	\chi\!\left(\frac{N(x)\cdot l}{Ra}\right).
	\]
	The matrix $K_a$ is positive semidefinite and
	\begin{equation}\label{eq:gramnorm}
		\|T_ac\|_2^2=\langle K_ac,c\rangle,\qquad
		\|T_a\|^2=\|T_a^*T_a\|=\|T_aT_a^*\|.
	\end{equation}
	\subsection{Kernel decay for signed and central bands}
	\begin{proposition}\label{prop:kernel}
		For every signed band $a_0\le a\le3$, and for the central band
		at $a=a_0$,
		\begin{equation}\label{eq:kernel}
			|K_a(k,l)|\lesssim_L a(1+a^2|k-l|)^{-L}
			\qquad(L=0,1,2,\ldots).
		\end{equation}
	\end{proposition}
	
	\begin{proof}
		We use the normal cutoffs and the equal-radius condition
		$|k|=|l|=R$ to keep the phase away from stationary points, and
		then integrate by parts. The signed bands and the central band
		require different geometric observations.
		
		\smallskip
		\noindent\emph{Amplitude bounds and the direct estimate.}
		Write $A=A_{a,k,l}$. On its support,
		$|x|\le2$ and $|x\cdot k'-k_4|\le4Ra$, where $k=(k',k_4)$.
		If $a\le a_*$, then
		\[
		1=\frac{|k|}{R}
		\le\frac{|k'|+|k_4|}{R}
		\le3\frac{|k'|}{R}+4a,
		\qquad |k'|\ge R/6.
		\]
		Thus the normal cutoff confines $x$ to a slab of thickness at
		most $8Ra/|k'|\lesssim a$. Its intersection with $B(0,2)$ has
		measure $O(a)$. For $a\ge a_*$, the fixed support of $\eta$
		gives the same bound with a larger fixed constant. Since $N$
		is affine and $k/R,l/R$ are unit vectors, differentiating a band
		cutoff costs at most $a^{-1}$ per derivative. Consequently,
		\begin{equation}\label{eq:amplitude-bounds}
			|\operatorname{supp}A|\lesssim a,\qquad
			\|\partial^\beta A\|_\infty\lesssim_\beta a^{-|\beta|}.
		\end{equation}
		In particular, $|K_a(k,l)|\lesssim a$. This proves the assertion
		when $k=l$ or $da^2\le8$, where $d=|l-k|$.
		
		For the remaining range $da^2>8$, normalize the frequency
		difference and the phase by setting
		\[
		\begin{gathered}
			v=\frac{l-k}{d}=(v',v_4),\qquad
			\phi(x)=v'\cdot x+\frac{v_4}{2}|x|^2,\\
			q(x)=\nabla\phi(x)=v'+v_4x,\qquad
			K_a(k,l)=\int_{\R^3}A(x)\ee(d\phi(x))\dd x.
		\end{gathered}
		\]
		Let $\nu(x)=N(x)/|N(x)|$. The identity
		$v=(q,0)-v_4N$ shows that stationarity, $q=0$, is equivalent
		to $l-k$ being parallel to the normal $N(x)$.
		Projecting the same identity onto $\nu^\perp$ gives
		\begin{equation}\label{eq:gradient}
			|q(x)|\ge
			|\operatorname{proj}_{\nu^\perp}v|
			=\sqrt{1-(v\cdot\nu)^2}.
		\end{equation}
		
		\smallskip
		\noindent\emph{Signed bands: the common sign excludes stationarity.}
		On the support of $A$, the two numbers $N(x)\cdot k$ and
		$N(x)\cdot l$ have the same sign and absolute value at least $Ra$.
		Hence
		\[
		|N(x)\cdot(k+l)|
		=|N(x)\cdot k|+|N(x)\cdot l|\ge2Ra.
		\]
		On the other hand, the equal-radius condition gives
		\[
		(l-k)\cdot(k+l)=|l|^2-|k|^2=0.
		\]
		A stationary point would therefore force $N(x)\cdot(k+l)=0$,
		a contradiction. This is the reason for separating the positive
		and negative normal bands.
		
		To quantify this observation, put $m=(k+l)/|k+l|$, which is
		well defined whenever the amplitude support is nonempty.
		Then $m\perp v$, and $|k+l|\le2R$, $|N(x)|\le\sqrt5$ give
		\[
		|\nu\cdot m|
		=\frac{|N(x)\cdot(k+l)|}{|N(x)|\,|k+l|}
		\ge\frac{a}{\sqrt5}.
		\]
		Because $m,v$ are orthonormal, $|\nu\cdot m|^2+|\nu\cdot v|^2\le1$.
		Combining this with \eqref{eq:gradient} yields
		\[
		|q(x)|\ge|\nu\cdot m|\ge a/3
		\qquad\hbox{on }\operatorname{supp}A.
		\]
		
		\medskip
		\noindent\emph{The central band.} Write  $D_0=R^{2/3}$.
		Here the signs may differ, but $a=a_0$ and the cutoff gives
		$|N(x)\cdot k|,|N(x)\cdot l|\le2Ra$.
		Since $|N(x)|\ge1$, subtraction gives
		\[
		|\nu\cdot v|
		\le\frac{|N(x)\cdot(l-k)|}{d}
		\le\frac{4Ra}{d}
		=\frac{4}{da^2}<\frac12.
		\]
		The equality uses the scale identity $Ra_0=D_0=a_0^{-2}$;
		the last inequality is precisely the range $d>8D_0$ now under
		consideration. Thus \eqref{eq:gradient} gives $|q|\ge\sqrt3/2$.
		The central cutoff can allow stationary pairs only at separations
		$O(D_0)$, where the direct estimate already suffices.
		
		\Needspace{6\baselineskip}
		\smallskip
		\noindent\emph{Integration by parts and the scale $da^2$.}
		Since $a\le3$, both cases give $|q|\ge a/4$ on the amplitude support.
		Let $V$ agree with $q/|q|^2$ on a neighborhood of that support
		and extend it smoothly by cutting it off where $|q|\le a/8$.
		Since $q$ is affine with bounded derivatives,
		\[
		\|\partial^\beta V\|_\infty\lesssim_\beta a^{-1-|\beta|},
		\qquad
		\frac{1}{2\pi i d}V\cdot\nabla\ee(d\phi)=\ee(d\phi)
		\]
		near $\operatorname{supp}A$. Set
		$\mathcal Df=\operatorname{div}(Vf)
		=V\cdot\nabla f+(\operatorname{div}V)f$.
		Each integration by parts contributes $d^{-1}$ and replaces
		the amplitude by $\mathcal D$ applied to it. The inverse
		gradient costs $a^{-1}$, and differentiation costs one more
		factor $a^{-1}$, explaining the scale $da^2$.
		
		More precisely, the product rule, the bounds for $V$, and
		\eqref{eq:amplitude-bounds} give inductively
		\[
		\|\partial^\beta\mathcal D^j A\|_\infty
		\lesssim_{j,\beta}a^{-2j-|\beta|},
		\qquad j=0,1,2,\ldots.
		\]
		These functions remain supported in $\operatorname{supp}A$.
		There are no boundary terms, since $A$ is compactly supported.
		After $L$ integrations by parts,
		\[
		K_a(k,l)=\left(-\frac1{2\pi i d}\right)^L
		\int_{\R^3}\mathcal D^L A(x)\ee(d\phi(x))\dd x,
		\]
		and hence
		\[
		|K_a(k,l)|
		\lesssim_L d^{-L}\|\mathcal D^L A\|_1
		\lesssim_L d^{-L}a^{1-2L}
		=a(da^2)^{-L}.
		\]
		Combining this with the direct estimate for $da^2\le8$
		proves \eqref{eq:kernel}.
	\end{proof}
	
	\subsection{The scalar and block Schur tests}
	\begin{lemma}[Schur test]\label{lem:Schur}
		Let $K_{ij}:H_j\to H_i$ be operators between Hilbert spaces.
		If
		\[
		\sup_i\sum_j\|K_{ij}\|\le A,\qquad
		\sup_j\sum_i\|K_{ij}\|\le B,
		\]
		then $K((v_i)_i)=(\sum_jK_{ij}v_j)_i$ has norm at most
		$\sqrt{AB}$ on the Hilbert direct sum
		$\bigoplus_iH_i$, whose squared norm is $\sum_i\|v_i\|^2$.
		For scalar matrices take $H_i=\mathbb C$ and replace norms by
		absolute values.
	\end{lemma}
	\begin{proof}
		For vectors $v=(v_j)$ and $w=(w_i)$, Cauchy--Schwarz gives
		\[
		\begin{aligned}
			|\langle Kv,w\rangle|
			&\le\sum_{i,j}\|K_{ij}\|\,\|v_j\|\,\|w_i\|\\
			&\le\left(\sum_{i,j}\|K_{ij}\|\|v_j\|^2\right)^{1/2}
			\left(\sum_{i,j}\|K_{ij}\|\|w_i\|^2\right)^{1/2}
			\le\sqrt{AB}\,\|v\|\,\|w\|.
		\end{aligned}
		\]
	\end{proof}
	
	\subsection{Combining frequency caps}
	Keep $\kappa=1/32$. For a small band scale $a$, partition $\R^4$ into
	half-open cubes $Q_i=D(i+[0,1)^4)$, $i\in\Z^4$.
	Every occupied cube lies in a ball $B(Re_i,4D)$, by choosing
	$Re_i\in\E_R\cap Q_i$. Put
	$H_i=\ell^2(\E_R\cap Q_i)$ and
	$K_{ij}=P_{Q_i}T_a^*T_aP_{Q_j}:H_j\to H_i$.
	Since $\sum_jP_{Q_j}=I$, these blocks assemble the Gram operator:
	\[
	\begin{aligned}
		T_a&=\sum_jT_aP_{Q_j},\\
		T_a^*T_a&=\sum_{i,j}(T_aP_{Q_i})^*(T_aP_{Q_j})
		=\sum_{i,j}K_{ij}.
	\end{aligned}
	\]
	Equivalently, with $c_j=P_{Q_j}c$,
	\[
	\|T_ac\|_2^2=\sum_{i,j}\langle K_{ij}c_j,c_i\rangle.
	\]
	Since
	$K_{ij}=(T_aP_{Q_i})^*(T_aP_{Q_j})$, Proposition~\ref{prop:local}
	gives
	\begin{equation}\label{eq:blocklocal}
		\|K_{ij}\|\le\|T_aP_{Q_i}\|\,\|T_aP_{Q_j}\|
		\lesssim A_D:=D^{-\kappa}.
	\end{equation}
	
	For $k\in Q_i$ and $l\in Q_j$, the offsets
	$k/D-i,l/D-j\in[0,1)^4$ have distance at most $2$, so
	\[
	|i-j|\le D^{-1}|k-l|+2=a^2|k-l|+2.
	\]
	Proposition~\ref{prop:kernel} therefore bounds every Gram matrix entry by
	\[
	|K_a(k,l)|\lesssim_L a(1+|i-j|)^{-L}.
	\]
	For $v\in H_j$, Cauchy--Schwarz in the inner sum gives
	\[
	\begin{aligned}
		\|K_{ij}v\|_2^2
		&=\sum_{k\in\E_R\cap Q_i}
		\left|\sum_{l\in\E_R\cap Q_j}K_a(k,l)v_l\right|^2\\
		&\le\left(\sum_{k\in\E_R\cap Q_i}
		\sum_{l\in\E_R\cap Q_j}|K_a(k,l)|^2\right)\|v\|_2^2\\
		&\lesssim_L a^2D^6(1+|i-j|)^{-2L}\|v\|_2^2.
	\end{aligned}
	\]
	Taking square roots yields
	\begin{equation}\label{eq:blockfar}
		\|K_{ij}\|\lesssim_L aD^3(1+|i-j|)^{-L}.
	\end{equation}
	This also covers neighboring cubes. Since $a=D^{-1/2}\le A_D$,
	the two estimates give
	\[
	\|K_{ij}\|\lesssim_L A_D\min\{1,D^3(1+|i-j|)^{-L}\}.
	\]
	Use $\min(1,t)\le t^\theta$ with $\theta=\kappa/8$ and
	choose $L$ so that $L\theta\ge5$. Since
	$-\kappa+3\theta=-5\kappa/8\le-\kappa/2$, we obtain
	\[
	\|K_{ij}\|\lesssim D^{-\kappa/2}(1+|i-j|)^{-5}.
	\]
	Since $\sum_{n\in\Z^4}(1+|n|)^{-5}<\infty$, the block norms satisfy
	\[
	\sup_i\sum_j\|K_{ij}\|+\sup_j\sum_i\|K_{ij}\|
	\lesssim D^{-\kappa/2}.
	\]
	Applying the block Schur test, Lemma~\ref{lem:Schur}, to $(K_{ij})$
	on $\ell^2(\E_R)=\bigoplus_iH_i$ yields
	\begin{equation}\label{eq:globalband}
		\|T_a\|^2=\|T_a^*T_a\|\lesssim D^{-\kappa/2}.
	\end{equation}
	
	\subsection{Summing the bands}\label{sec:bandsum}
	We now sum the decomposition from Section~\ref{sec:band-definition}.
	By \eqref{eq:globalband},
	\[
	\|T_{a,\pm}\|\lesssim D^{-\kappa/4}=a^{1/64},\qquad
	\|T_0\|\lesssim D_0^{-\kappa/4}=R^{-1/192}.
	\]
	The dyadic sum for $a\le a_*$ is a convergent geometric series,
	uniformly in $R$. There are only boundedly many bands with
	$a_*<a\le3$. For these, Proposition~\ref{prop:kernel} gives
	$|K_a(k,l)|\lesssim_{L,a_*}(1+|k-l|)^{-L}$. Choosing $L>4$,
	unit-ball packing bounds the scalar row and column sums:
	\[
	\sup_k\sum_l|K_a(k,l)|+\sup_l\sum_k|K_a(k,l)|
	\lesssim\sum_{h\in\Z^4}(1+|h|)^{-L}<\infty.
	\]
	The scalar Schur test (Lemma~\ref{lem:Schur} with $H_i=\mathbb C$)
	therefore gives $\|T_a\|^2=\|K_a\|\lesssim1$.
	The triangle inequality proves $\|Tc\|_2\lesssim\|c\|_2$.
	Since $\eta=1$ on $[-1,1]^3$, Parseval in \eqref{eq:eigenfunction}
	now gives Theorem~\ref{thm:main}.

	\section*{Use of artificial intelligence}
	The author's interest in this problem dates back to 2020,
	when he explored an approach based on the Bourgain--Rudnick machinery
	and the four-dimensional cap-energy estimate proved in joint work
	with Bourgain \cite[Section~4]{BD}. While that approach did not yield the
	desired restriction bound, it subsequently served as a starting point
	for an extended interaction with ChatGPT, in which the author guided
	the development of increasingly efficient arguments. This process
	led to a proof slightly more complicated than the one presented here.
	The author further simplified the argument and wrote the present
	version.


	\label{page:last}
\end{document}